\documentclass[11pt]{amsart}
\usepackage[margin=1in]{geometry}
\usepackage[T1]{fontenc}
\usepackage{lmodern}
\usepackage[indent=0pt]{parskip}

\usepackage{tikz}
\usetikzlibrary{calc}
\usetikzlibrary{arrows.meta}

\usepackage{xcolor}
\usepackage{graphicx}
\usepackage{amsmath,amssymb,mathtools,microtype}
\usepackage[hidelinks]{hyperref}
\newtheorem{thm}{Theorem}[section]
\newtheorem{lem}[thm]{Lemma}
\newtheorem{prop}[thm]{Proposition}
\newtheorem*{KMthm}{Kronheimer--Mrowka's theorem}

\makeatletter
\newcommand{\namedlabel}[2]{%
  \begingroup
  \def\@currentlabel{#2}%
  \label{#1}%
  \endgroup
}
\makeatother

\theoremstyle{remark}
\newtheorem{rem}[thm]{Remark}
\newcommand{\C}{\mathbb C}
\newcommand{\R}{\mathbb R}
\newcommand{\CP}{\mathbb{CP}}
\DeclareMathOperator{\rank}{rank}
\title[polyanalytic]{On the number of zeros of polyanalytic polynomials}
\author{Erik Lundberg}
\address{Department of Mathematical Sciences, Florida Atlantic University,
Boca Raton, FL 33431, USA}
\email{elundber@fau.edu}
\author{Zvi Rosen}
\address{Department of Mathematical Sciences, Florida Atlantic University,
Boca Raton, FL 33431, USA}
\email{rosenz@fau.edu}
\author{Teng Zhang}
\address{School of Mathematics and Statistics, Xi'an Jiaotong University,
Xi'an 710049, P.~R.~China}
\email{teng.zhang@stu.xjtu.edu.cn}
\date{October 16, 2026}
\begin{document}
\begin{abstract}
We prove that a complex polynomial $f(z,\bar z)$ of degree $n$ in $z$ and degree $m<n$
in $\bar z$, with finitely many zeros, has at most $(4m+1)n-5m$ distinct
zeros in $\C$. 
The proof uses an approach suggested by GPT-6.1 based on a topological result of Kronheimer and Mrowka (Thom's Conjecture).
We provide examples showing that the bound is asymptotically
sharp as $n\to\infty$ for each fixed $m$.
\end{abstract}
\maketitle

\section{Introduction}

In this paper we study the number of zeros of polynomials in $z$ and $\bar z$.

For a bivariate polynomial $f\in\C[z,w]$, set $F_f(z)=f(z,\bar z)$ and let
$N(f)=\#F_f^{-1}(0)$ denote the number of distinct zeros. Let $n=\deg_z f$ and $m=\deg_w f$, i.e., $(n,m)$ is the bidegree of $f$.
The problem of bounding $N(f)$ goes back to Sheil-Small \cite[Section 1.2]{SS02}, \cite{SS92}.
The particular cases of harmonic polynomials
\cite{BHS95,Wil98,KS03,Gey08,HLLM15,LLL15,KLS18,Lun23}
and logharmonic polynomials
\cite{BH,AH01,BL10,KLP24,LL}
have received a lot of attention.
Khavinson, Lundberg, and Perry
\cite[Problem~1.9 and Proposition~1.10]{KLP24}
formulated the general problem of bounding $N(f)$ in terms of $n,m$, and they used resultants to prove the bound $N(f) \leq n^2+m^2$
when $f$ and $f^*(z,w)=\overline{f(\bar w,\bar z)}$ are coprime.

A central question is whether, for fixed $m$, this quadratic
dependence on $n$ can be replaced by a linear bound.
Eremenko \cite[p.~1]{Ere26} stated the following problem,
generalizing several open questions posed for particular subclasses of
polyanalytic polynomials \cite{Wil98,BH,LLL15,Oka18,SZ24}.

{\bf Problem.}
{\it For each \(m\ge1\), do there exist constants \(A_m,B_m\) such that
\(N(f)\le A_m n+B_m\) for every polynomial of bidegree \((n,m)\),
\(n>m\), with finite zero set?}

We prove the following result giving an affirmative answer.

\begin{thm}\label{thm:linear}
Let $n>m\ge1$ and let $f\in\C[z,w]$ with
$\deg_z f = n$ and $\deg_w f  = m$. If $N(f)<\infty$, then
\begin{equation}\label{eq:mainbound}
 N(f)\le(4m+1)n-5m.
\end{equation}
\end{thm}
Part of our proof (see Section \ref{sec:regular}) also shows that the above-mentioned bound of \cite{KLP24} holds under the finiteness assumption $N(f) < \infty$, without the additional coprimality assumption. Hence, under the assumptions in Theorem \ref{thm:linear} we actually have the bound $N(f)\le\min\{n^2+m^2,(4m+1)n-5m\}$.

In the case $m=1$, our theorem recovers Khavinson and Neumann's $5n-5$ bound \cite{KN06} for the valence of certain rational harmonic mappings.
More precisely, their result bounds the number of solutions of equations of the form
\[
 \frac{p(z)}{q(z)}-\bar z = 0,
\]
where $p,q$ are coprime univariate complex polynomials, and $\max\{\deg p,\deg q\}=n\ge2$.
Multiplying by $q(z)$ relates this to the zero set of a polyanalytic polynomial
$p(z)-\bar zq(z)$ of bidegree $(n,1)$.
The $5n-5$ bound was conjectured by astronomer Rhie \cite{Rh01,Rh03} in connection with image-counting problems in gravitational lensing (she also showed that $5n-5$ is sharp using constructions motivated by astrophysical considerations).
The proof in \cite{KN06} makes essential use of tools from holomorphic dynamics.

The bound \eqref{eq:mainbound} holds for the special classes of harmonic polynomials $p(z)+\overline{q(z)}$ with
$\deg p=n>\deg q=m$, as well as logharmonic polynomials $p(z) \overline{q(z)}$. For harmonic polynomials, Wilmshurst proved the upper bound of $n^2$ for the valence using B\'ezout's Theorem \cite{Wil98}, and he conjectured the refined (linear in $n$) bound $m(m-1)+3n-2$. This was confirmed in the case $m=1$ by Khavinson and \'Swi\k{a}tek \cite{KS03}.
Counterexamples for $m=n-3$ large were provided by Lee, Lerario, and Lundberg in \cite{LLL15}, where they conjectured the weakened bound of $2m(n-1) + n$.
For zero-counting results in terms of total degree and a
comparison between harmonic and polyanalytic valence,
see S\`ete and Zur \cite{SZ24}.

For generic coefficients, the proof of Theorem \ref{thm:linear} constructs a smooth embedded oriented
surface $\Sigma\subset\CP^2$ with
\begin{equation}\label{eq:goal}
 g(\Sigma)=(n-1)(m-1)+\frac{n^2+m^2-N}{2},
 \qquad d:=\deg\Sigma=n-m,
\end{equation}
where $N=N(f)$. The Kronheimer--Mrowka theorem (formerly, the Thom Conjecture) stated below gives a lower bound for
its genus, hence an upper bound for $N$.

\begin{KMthm}[{\cite[Theorem~1]{KM94}}] \namedlabel{thm:thom}{Kronheimer--Mrowka theorem}

If $\Sigma\subset\CP^2$ is a smooth closed connected oriented embedded surface
of degree $d>0$, then
\[
 g(\Sigma)\ge\frac{(d-1)(d-2)}2,
\]
where the degree $d$ can be defined as the signed intersection number
with a transverse complex projective line:
\[
 d=\sum_{p\in\Sigma\cap\ell}\operatorname{sign}(p).
\]
\end{KMthm}

The sign at an intersection point $p$ is positive when the two oriented tangent
planes give the ambient orientation, and negative otherwise.

\begin{rem}\label{rem:homology}
In \cite{KM94}, the degree $d$ is defined in terms of homology. Recall that
$H_2(\CP^2;\mathbb Z)=\mathbb Z[\ell]$, where $\ell$ is a complex
line \cite[\S2.2, Theorem~2.35 and p.~140]{Hat02}. The homological degree $d$ of a surface $S$ is defined by $[S]=d[\ell]$.  This coincides with the definition using signed intersection, since one has $S\cdot\ell=d$ because
$\ell\cdot\ell=1$ \cite[Prop. 1.2.5]{GS99}.
\end{rem}

We begin by explaining the construction of the surface $\Sigma$ for the case of generic coefficients in Section \ref{sec:generic}. Then, in Section \ref{sec:regular}, we reduce the general case to the generic case by a degree-preserving perturbation that regularizes degenerate zeros without decreasing the cardinality of the zero set.
Section \ref{sec:examples} provides examples showing that
Theorem \ref{thm:linear} is asymptotically sharp as $n\to\infty$
for each fixed $m$ (note: these examples are not harmonic, so they do not violate the above mentioned $2m(n-1)+n$ conjecture).

\noindent\textbf{Acknowledgments and AI tools disclosure.} 
The idea to use a topological approach based on the theorem of Kronheimer and Mrowka was suggested in an interaction with GPT-6.1. The model was also used in subsequent discussions to explore and refine the argument and to help create Figure \ref{fig:psi-cartoon}. The final proofs were written and checked by the authors.  This work was supported in part by a grant of access to OpenAI
models through the ChatGPT for Academic Researchers program.
The first-named author also acknowledges support from Simons Foundation grant 712397. 


\section{Proof of Theorem \ref{thm:linear} in the generic case}\label{sec:generic}

Bihomogenize $f$ to bidegree $(n,m)$, and let $C\subset \CP^1_z\times\CP^1_w$
be its zero curve.
Throughout this section, we assume that $C$ is smooth
and avoids $(\infty,\infty)$.
We are interested in the zeros of $f$ on the affine chart that lie on the non-holomorphic subspace
\[
 L=\{(z,\bar z):z\in\CP^1\},
\]
which is a smooth real two-sphere. The boundary of the affine chart is
the union of $z=\infty$ and $w=\infty$; it meets $L$ only at
$(\infty,\infty)$. Since $C$ avoids this point, we have
\[
 N=N(f)=\#(C\cap L).
\]

Define a reflection (antiholomorphic involution on $\CP^1\times\CP^1$) by
\[
 \tau(z,w)=(\bar w,\bar z).
\]
Note that $L$ is precisely the fixed set of $\tau$. The reflected curve $\tau C=\{\tau(x):x\in C\}$ is the zero curve
of the bihomogenization of $f^*$ to bidegree $(m,n)$, where
$f^*(z,w)=\overline{f(\bar w,\bar z)}$.
Since $\tau$ fixes $L$ pointwise, every point of $C\cap L$ also
lies on $\tau C$, so we have
\begin{equation}
    C\cap L\subset C\cap\tau C.
\end{equation}
Assume also that $C$ meets both $L$ and $\tau C$ transversely.

The following lemma shows that proving 
Theorem~\ref{thm:linear} reduces to establishing the bound under the generic assumptions specified above (which are summarized in the lemma statement).  To clarify, we recall that zeros are not counted with multiplicity.

\begin{lem}\label{lem:generic-reduction}
The above transversality and genericity conditions, namely, that $C$ is smooth, avoids $(\infty,\infty)$, and meets both $L$ and $\tau C$
transversely, hold for an open dense set of coefficients.
Moreover, every polynomial $f$ with $N(f)<\infty$ can be perturbed to a polynomial of the same bidegree satisfying these genericity conditions and having at least $N(f)$ distinct zeros.
\end{lem}

The proof of the lemma is deferred to Section~\ref{sec:regular} in favor of first presenting the proof of Theorem \ref{thm:linear} in the generic case.
Hence, the genericity conditions stated in Lemma \ref{lem:generic} are assumed throughout this section.

\subsection{Connectedness of \texorpdfstring{$C$}{C}, genus of \texorpdfstring{$C$}{C}, and intersection number of \texorpdfstring{$C$}{C} and \texorpdfstring{$\tau C$}{tau C}}\label{sec:curve}

In order to compute the genus of the surface $\Sigma$ constructed in the subsequent sections, we will need a few basic results from algebraic geometry:

\begin{lem}[{\cite[Chapter~V, Example~1.4.3]{Har77}}]
\label{lem:intersection}
Two curves in $\CP^1\times\CP^1$ of bidegrees $(n_1,m_1)$
and $(n_2,m_2)$, with no common component, have intersection
number $n_1m_2+m_1n_2$.
\end{lem}

Applying Lemma~\ref{lem:intersection} to $C$ and $\tau C$, while noting (by transversality) that every intersection multiplicity is one
\cite[Chapter~V, Theorem~1.1(1)]{Har77}, gives
\begin{equation}\label{eq:intersection}
 \#(C\cap\tau C)=n^2+m^2.
\end{equation}

\begin{lem}\label{lem:connected}
Every curve $C\subset \CP^1\times \CP^1$ of
bidegree $(n,m)$ with $n,m>0$ is connected.
\end{lem}

\begin{proof}
Each irreducible component of $C$ is connected in the usual
topology \cite[Ch.~7, \S2.2, Lemma~7.3]{Shafarevich}.
Thus, if $C$ were disconnected, its irreducible components
could be grouped into two nonempty disjoint algebraic
curves $A$ and $B$. Their bidegrees $(a,b)$ and $(c,d)$,
counting multiplicities, have nonnegative entries and satisfy
\[
a+c=n,\qquad b+d=m.
\]
Since $A$ and $B$ are disjoint, the intersection formula provided in Lemma \ref{lem:intersection} gives
\[
ad+bc=0.
\]
Neither bidegree is $(0,0)$, so this forces either
$a=c=0$ or $b=d=0$, contradicting $n,m>0$.
\end{proof}

\begin{lem}[{\cite[Chapter~V, Example~1.5.2]{Har77}}]
\label{lem:genus}
A smooth connected curve of bidegree $(n,m)$ in
$\CP^1\times\CP^1$ has genus $(n-1)(m-1)$.
\end{lem}

Applying Lemma~\ref{lem:genus} to $C$ gives
\begin{equation}\label{eq:genus}
 g(C)=(n-1)(m-1).
\end{equation}

\subsection{Mapping \texorpdfstring{$C$}{C} to an immersed surface in the projective plane}\label{sec:immersion}

The map $(z,w) \mapsto (w+\bar z,w\bar z)$  on the affine charts extends to a smooth map $\Psi: \CP^1\times\CP^1\to\CP^2$ by 
\begin{equation}\label{eq:map}
 \Psi([Z_0:Z_1],[W_0:W_1])
 =[\bar Z_0W_0:\bar Z_0W_1+\bar Z_1W_0:\bar Z_1W_1].
\end{equation}
The coordinates never vanish simultaneously. 
They are the coefficients of the quadratic
$(\bar Z_0T+\bar Z_1S)(W_0T+W_1S)$.
Unique factorization shows that the fibers are exactly the
$\tau$-orbits, including at infinity. In particular, $\Psi\circ\tau=\Psi$.

\begin{figure}[htbp]
\centering
\includegraphics[width=\textwidth]{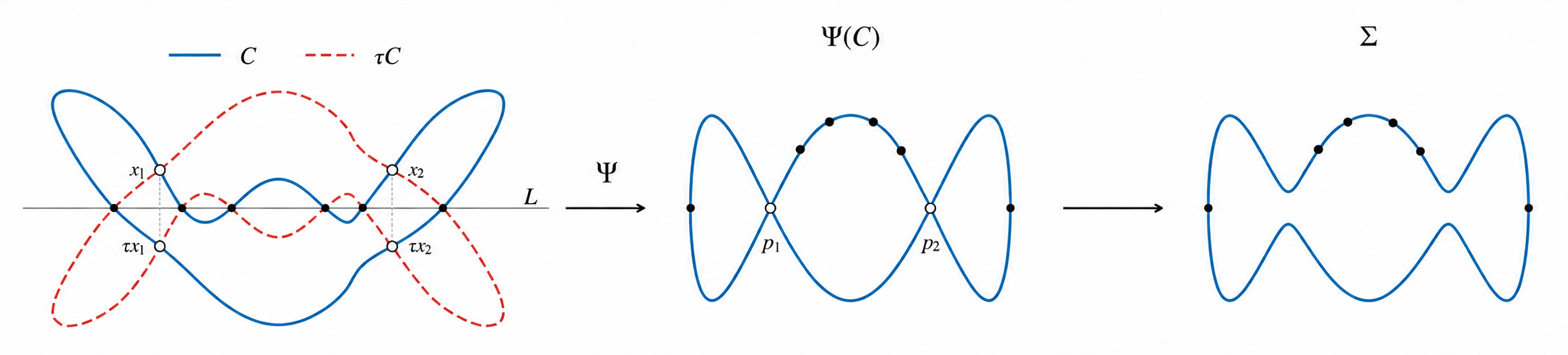}
\caption{Schematic of the construction of the surface $\Sigma$ (with real surfaces depicted as curves).
The six intersections of $C \cap \tau C$ that lie on $L$ are mapped to six smooth points on $\Psi(C)$. The two pairs $(x_j,\tau x_j)$ off $L$ are mapped to the two double points $p_j=\Psi(x_j)=\Psi(\tau x_j)$,
which are then resolved to obtain $\Sigma$.}
\label{fig:psi-cartoon}
\end{figure}

Figure~\ref{fig:psi-cartoon} illustrates this identification and the
subsequent smoothing of the double points.

\begin{lem}\label{lem:immersion}
The restriction $\Psi|_C$ is an immersion, and the image $\Psi(C)$ has precisely
\begin{equation}\label{eq:delta}
 \delta=\frac{n^2+m^2-N}{2}
\end{equation}
double points and no other multiple points.
Moreover, every double point is transverse and has positive intersection sign.
\end{lem}
\begin{proof}
In finite affine coordinates, its real differential is
\begin{equation}\label{eq:differential}
 d\Psi_{(z,w)}(a,b)=(b+\bar a,\bar z\,b+w\bar a).
\end{equation}
It is invertible off $L$, because when both components are zero, $b + \bar{a} = 0$ implies $\bar{z} b + w\bar{a} = (\bar{z}-w)b = 0$. 
On $L$,
\[
 \ker d\Psi=\{(a,-\bar a):a\in\C\},\qquad
 TL=\{(a,\bar a):a\in\C\}.
\]
If $v\in TC \cap\ker d\Psi$, then $iv \in TC  \cap  TL$, since $TC$
is a complex line and  \[iv = (ia,-i\bar{a}) = (ia, \overline{ia}) \in TL.\] Transversality of $C$ and $L$ implies $v=0$, which proves injectivity of the differential throughout the affine chart.

 Near $(\infty,w)$, use the coordinate $t=1/z$, i.e., $Z_0 = t, Z_1 = 1$; then, the image of $\Psi$ is $[\bar{t}:1+ \bar{t}w : w]$. Normalize with respect to the second coordinate, since it is nonzero near $t=0$,
 to write
\[
 \Psi(t,w)=\left(\frac{\bar t}{1+\bar t w},
                  \frac{w}{1+\bar t w}\right).
\]
Its differential at $t=0$ is $(a,b)\mapsto(\bar a,b-w^2\bar a)$,
which is invertible. The case $(z,\infty)$ follows from
$\Psi\circ\tau=\Psi$, and $C$ avoids $(\infty,\infty)$.
Thus $\Psi|_C$ is an immersion everywhere. The local normal form for
immersions \cite[Chapter~1, \S3]{GP74} and compactness of $C$ imply
that its image $\Psi(C)$ is smooth at every point having a unique preimage.

As observed above, unique factorization implies that fibers of $\Psi$ are precisely $\tau$-orbits, which are of size 1 or 2. Therefore, distinct points $x \neq y$ in $C$ share an image exactly when $y = \tau x$
with $x\in(C\cap\tau C)\setminus L$. The $\delta$ double points of the curve $\Psi(C)$ are
in one-to-one correspondence with (unordered) pairs $\{x,y\}$ of this type. Hence, grouping the points in
$C \cap \tau C$ by fiber of $\Psi$, we have $n^2 + m^2 = 2\delta +N$, which implies \eqref{eq:delta}.

Fixing one of these double points, we now verify transversality and positivity of its sign (the
sign of a double point is defined by the usual intersection-sign convention 
\cite[Chapter~3, \S3, p.~113]{GP74} applied to the two local image branches).
This requires showing that the following procedure produces a positively oriented basis for $T_{\Psi(x)}\CP^2$: take positively oriented bases of
$T_xC$ and $T_{\tau x}C$, map them by $d\Psi_x$ and $d\Psi_{\tau x}$,
respectively, then concatenate the resulting pair.

Denote the chosen bases for $T_xC$ and $T_{\tau x}C$ by $\mathcal B_1$ and $\mathcal B_2$, respectively.
The identity $\Psi\circ\tau=\Psi$ and the chain rule at $\tau x$ give
\[
 d\Psi_{\tau x}=d(\Psi\circ\tau)_{\tau x}
 =d\Psi_x\circ d\tau_{\tau x}.
\]
In particular, the concatenation of image vectors can be expressed as
\[
 \bigl(d\Psi_x\mathcal B_1,d\Psi_{\tau x}\mathcal B_2\bigr)
 =d\Psi_x\bigl(\mathcal B_1,d\tau_{\tau x}\mathcal B_2\bigr).
\]
Note that $\left( \mathcal B_1,d\tau_{\tau x}\mathcal B_2 \right)$ is a basis of
$T_x(\CP^1\times\CP^1)$, since $C$ and $\tau C$ are transverse.  Moreover, it has negative orientation, since  $d \tau$ reverses the orientation of the complex line $T_{\tau x}C$ relative to the complex orientation of
$T_x(\tau C)$.
Meanwhile, $d\Psi_x$ has negative
real determinant: $\Psi$ is conjugation in the first factor followed
by a holomorphic local diffeomorphism. Thus, the orientation sign of the
image basis is $(-1)(-1)=+1$, proving positivity.
\end{proof}

\subsection{Resolving the double points of \texorpdfstring{$\Psi(C)$}{Psi(C)} to produce the smooth surface \texorpdfstring{$\Sigma$}{Sigma}}\label{sec:genus}
Using the positivity of the double points of $\Psi(C)$ stated in Lemma \ref{lem:immersion}, we can resolve each double point  by following the local surgery in \cite[\S2.1, p.~38 and Remark~2.1.2]{GS99}.
Namely, within a small ball centered at a given double point, $\Psi(C)$ consists of a pair of transverse disks.  Cut out the interiors of both disks and replace them by an annulus glued smoothly to the surface along the two circular boundary components where the disks were removed. As in \cite[\S2.1, p.~38 and Remark~2.1.2]{GS99}, the orientation of the surface (the orientation on $\Psi(C)$ is the one induced by the orientation on $C$) extends across the annulus. Denote by $\Sigma$ the smooth surface obtained by resolving all double points (choosing the balls in which the smoothing is performed to be pairwise disjoint). Then $\Sigma$ is closed, connected, smoothly embedded, and oriented by the orientation inherited from $C$.

\begin{lem}\label{lem:smoothedgenus}
The genus of $\Sigma$ is
\begin{equation}\label{eq:smoothedgenus}
g(\Sigma)=g(C)+\delta
 =(n-1)(m-1)+\frac{n^2+m^2-N}{2}.
\end{equation}
\end{lem}

\begin{proof}
Determining the genus of $\Sigma$ is done by inspecting the above local surgery to compute its Euler characteristic $\chi (\Sigma)$.

Each ball in which the smoothing is performed meets $\Psi(C)$
in exactly the two local disks through its double point.
Let $S$ be the surface with boundary that remains after the interiors of these disks
are removed (before gluing the replacement annuli). 
The restriction of $\Psi$ maps $C$ minus the interiors of
$2\delta$ pairwise disjoint disks diffeomorphically onto $S$ (recall $\delta$ is the number of double points).
In particular, $S$ is connected. Since each disk removed decreases Euler characteristic by one, and Euler characteristic is a diffeomorphism invariant, we conclude
\[
 \chi(S)=\chi(C)-2\delta.
\]
The surface $\Sigma$ is obtained from $S$ by attaching $\delta$
annuli along its $2\delta$ boundary circles. 
Both an annulus and a circle have Euler characteristic zero,
so the inclusion--exclusion formula for Euler characteristic
\cite[\S2.2, Exercise~21, p.~157]{Hat02} gives
\[
 \chi(\Sigma)=\chi(S)=\chi(C)-2\delta.
\]

Since $C$ and
$\Sigma$ are closed connected oriented surfaces, the identity
$\chi=2-2g$ yields
\[
 2-2g(\Sigma)=2-2g(C)-2\delta,
\]
and hence $g(\Sigma)=g(C)+\delta$. This, together with
\eqref{eq:genus} and \eqref{eq:delta}, proves \eqref{eq:smoothedgenus}.
\end{proof}

The balls in which the smoothing is performed may be chosen
arbitrarily small; below, we choose them disjoint from the line
used to compute the degree.

\subsection{Computing the degree of \texorpdfstring{$\Sigma$}{Sigma}}\label{sec:degree}

In coordinates $[A:B:D]$, put
\[
 H_a=\{D-aB+a^2A=0\},\qquad a\in\C.
\]
Its defining equation pulls back in the finite chart to
\begin{equation}\label{eq:linepullback}
 w\bar z-a(w+\bar z)+a^2=(w-a)(\bar z-a),
\end{equation}
and globally $\Psi^{-1}(H_a)=\{w=a\}\cup\{z=\bar a\}$.
Choose $a$ so these fibers meet $C$ transversely at finite points,
$(\bar a,a)\notin C$, and $H_a$ avoids the double points.
Only finitely many values are excluded: the critical values of the two
projections of $C$, the values contributed by its boundary points and
$C\cap L$, and the roots of the line equation at the double points.
Choose the balls where the smoothing operation occurred to be disjoint from $H_a$.

The fiber $w=a$ meets $C$ in $n$ points, and $z=\bar a$ meets it in
$m$ points, by Lemma~\ref{lem:intersection}
(or by the Fundamental Theorem of Algebra).
With the orientation induced from $C$, the first $n$ intersections
with $H_a$ have sign $+1$, since $w|_C$ is holomorphic, and the other
$m$ have sign $-1$, since $\bar z|_C$ is antiholomorphic.
Indeed, in \eqref{eq:linepullback} the other factor is nonzero and
multiplies the normal differential by a nonzero complex scalar,
which preserves orientation. Smoothing has changed none of these
intersections or their signs. Hence
\begin{equation}\label{eq:degree}
 d:=\deg\Sigma=n-m>0.
\end{equation}

\subsection{Applying the genus bound in the \ref{thm:thom} (Thom Conjecture)}

\begin{proof}[Completion of the proof of Theorem~\ref{thm:linear} in the generic case] Applying 
the~\ref{thm:thom} to the smoothed surface $\Sigma$, together with \eqref{eq:smoothedgenus} and
\eqref{eq:degree}, gives
\[
 (n-1)(m-1)+\frac{n^2+m^2-N}{2}
 \ge\frac{(n-m-1)(n-m-2)}2.
\]
Rearranging yields the desired bound \eqref{eq:mainbound}.
\end{proof}

\section{Reduction of the general case to the generic case}\label{sec:regular}

In this section, we prove Lemma \ref{lem:generic-reduction}, showing that in proving the bound on $N(f)$, it was enough to consider only those polynomials $f(z,\bar{z})$ which (1) have nondegenerate zeros, (2) whose bihomogenization is smooth, (3) avoiding $(\infty, \infty)$, and (4) transverse to both $L$ and $\tau C$. In general, the strategy is to show that for any $f(z,\bar{z})$ failing one of these conditions, there is a nearby perturbed polynomial $\tilde{f}(z,\bar{z})$ satisfying the condition and having at least as many zeros. Since $\tilde{f}$ must satisfy the bound by the arguments in section 2, the original $f$ must satisfy the bound as well.

Related reductions to regular zeros in the harmonic setting
appear in \cite[\S4]{KS03} and
\cite[Proposition~3 and \S5]{BE10}.

The real Jacobian of $F_f$ is
\[
 J_f(z)=|f_z(z,\bar z)|^2-|f_w(z,\bar z)|^2.
\]
A zero is nondegenerate when $J_f\ne0$; the differential has rank zero
exactly when $f_z=f_w=0$. Nondegenerate zeros persist under small
coefficient perturbations. 
The main difficulty is to construct a perturbation which regularizes degenerate zeros without decreasing their cardinality.
This is handled in Lemma \ref{lem:regularization} which partly relies on the following local result.

\begin{lem}\label{lem:rankone}
Let $F$ be a real-analytic map from a neighborhood of $0\in\R^2$
to $\R^2$, with an isolated zero at $0$ and $\rank DF(0)=1$.
Let $G$ be real analytic near $0$, with
$G(0)\notin\operatorname{im}DF(0)$.
For any neighborhood $U$ of $0$, there exists $t_0>0$ such that, for every
$0<t<t_0$, the sum of the numbers of distinct nondegenerate zeros of $F+tG$ and
$F-tG$ in $U$ is at least two.
\end{lem}

\begin{proof}
By making an invertible linear change of coordinates in the domain, we may assume that $F_x(0)\ne0$.
Since $DF(0)$ has rank one and
$G(0)\notin\operatorname{im}DF(0)$, the vectors
$F_x(0)$ and $G(0)$ are linearly independent.
Thus, the real-analytic implicit function theorem
provides functions $x=X(y), s=S(y)$ solving the equation
\begin{equation}\label{eq:IFT}
 F(x,y)+s\, G(x,y)=0
\end{equation}
for $x,s$ in terms of $y$, i.e.,
\[
 (x,y,s)=(X(y),y,S(y)),
 \qquad X(0)=S(0)=0,
\]
parametrizes the nearby solutions of \eqref{eq:IFT} using $y$ as the parameter.
The function $S$ has an isolated zero at $0$: otherwise,
the points $(X(y),y)$ with $S(y)=0$ would give zeros of
$F$ accumulating at $0$. By analyticity,
\[
 S(y)=ay^\ell+O(y^{\ell+1})
\]
for some integer $\ell\ge1$ and $a\ne0$.
If $\ell$ is odd, each equation $S(y)=t$ and $S(y)=-t$
has a solution near $0$ for every sufficiently small $t>0$.
If $\ell$ is even, the equation
$S(y)=\operatorname{sgn}(a)t$ has two such solutions,
one on either side of $0$.
The corresponding points $(X(y),y)$ lie in $U$ when $t$
is sufficiently small. Hence $F+tG$ and $F-tG$ together
have at least two distinct zeros, counted separately
for the two signs.

Finally, to see that the zeros obtained are nondegenerate, we note that $S'(y)\ne0$ for all sufficiently small $y\ne0$.
Writing $H(x,y,s)=F(x,y)+sG(x,y)$ and differentiating
$H(X(y),y,S(y))=0$ gives
\[
 \det(H_x,H_y)=-S'(y)\det(H_x,G)\ne0,
\]
since $H_x$ and $G$ remain linearly independent near the origin.
\end{proof}

We will need the following special case of the transversality theorem from
\cite[Chapter~2, \S3, Transversality Theorem]{GP74}.
In the proof of Lemma \ref{lem:regularization}, this is applied with $U$ being $\mathbb R^2$ with a finite set $X_0$ removed.

\begin{prop}[Parametric transversality]\label{prop:parametric-transversality}
For an open set $U \subset \mathbb R^n$, let
\[
\mathcal F: U \times\mathbb R^m\to\mathbb R^n
\]
be smooth, and suppose that $D\mathcal F$ is surjective at every point of
$\mathcal F^{-1}(0)$. For $\xi\in\mathbb R^m$, set
\[
\mathcal F_\xi(x)=\mathcal F(x,\xi).
\]
Then, for almost every $\xi\in\mathbb R^m$, every zero of
$\mathcal F_\xi$ is nondegenerate; equivalently,
\[
D_x\mathcal F(x,\xi):\mathbb R^n\to\mathbb R^n
\]
is invertible whenever $\mathcal F_\xi(x)=0$.
\end{prop}

\begin{proof}
Apply the transversality theorem of \cite[Ch.~2, Sec.~3]{GP74}
with $X=U$, $S=\mathbb R^m$, $Y=\mathbb R^n$,
and $Z=\{0\}$. Since $T_0Z=\{0\}$, transversality of
$\mathcal F$ to $Z$ means exactly that $D\mathcal F$ is
surjective at every point of $\mathcal F^{-1}(0)$, as assumed.
The boundary condition is vacuous since $X$ has no boundary.
The theorem therefore implies that, for almost every $\xi \in S$,
$\mathcal F_\xi$ is transverse to $Z$; that is,
$D_x\mathcal F(x,\xi)$ is surjective at every zero of
$\mathcal F_\xi$. Being a linear map from $\mathbb R^n$
to itself, it is then invertible.
\end{proof}

\begin{lem}\label{lem:regularization}
Let \(f\in\C[z,w]\) have finitely many zeros on \(w=\bar z\).
Every coefficient neighborhood of \(f\) contains a polynomial
\(\widetilde f\) such that
\(\deg_z\widetilde f = \deg_z f\),
\(\deg_w\widetilde f = \deg_w f\),
and \(N(\widetilde f)\ge N(f)\), with all zeros nondegenerate.
\end{lem}

For a zero $\zeta$ of $f(z,\bar z)$, its \emph{order} is the least degree of a nonzero homogeneous term in the Taylor expansion of $f(z,w)$ at $(\zeta,\bar\zeta)$. Thus a zero of rank zero has order at least two. A zero of order one has nonzero differential, and hence is either rank one or full rank as a real map $\R^2\to\R^2$.

\begin{proof}
Let $X$ be the set of zeros of $f(z,\bar{z})$. Then, 
$X = X_0 \: \sqcup \: X_1 \: \sqcup \: X_2$
where $X_0:=$ the set of rank-0 zeros, $X_1 :=$ the set of rank-1 zeros, and $X_2 :=$ the set of nondegenerate zeros.

We will construct an iterative perturbation $\tilde{f}$ of $f$ with the same degree which (using $\tilde X = \tilde X_0 \sqcup \tilde X_1 \sqcup \tilde X_2$ to denote the corresponding partition of the zero set of $\tilde f$)
\begin{enumerate}
\item[(i)] strictly decreases the order of vanishing of zeros in $X_0$, possibly placing them in $\tilde X_1$ or $ \tilde X_2$;
\item[(ii)] replaces the zeros in $X_1$ with a nearby set of zeros in $\tilde X_2$ of equal or greater size; and
\item[(iii)] does not lose any zeros from $X_2$ but merely moves them to a nearby set of zeros in $\tilde X_2$ of the same size.
\item[(iv)] Moreover, $\tilde f$ has finitely many zeros, and all zeros of $\tilde f$ outside $X_0$ are nondegenerate, i.e., we have $\tilde X_0 \cup \tilde X_1 \subset X_0$.
\end{enumerate}

Note that this does not decrease the total number of zeros.
By iterating this perturbation, eventually $\tilde X_0$ and $ \tilde X_1$ both become empty, so that $\tilde{f}$ has only nondegenerate zeros as desired.  Hence, the lemma is proved as soon as we establish 
the following claim which furnishes the desired perturbation. 


Claim. {\em 
There exist arbitrarily small nonzero $\xi = (a,b)\in\C^2$ such that each of the maps
\[
f_{\pm \xi}(z,\bar z)
=f(z,\bar z)\pm\bigl(a f_z(z,\bar z)+b f_w(z,\bar z)\bigr),
\qquad z\in\C\cong\R^2,
\]
satisfies the properties (i), (iii), and (iv), and at least one of the two maps $f_\xi$ or $f_{-\xi}$ satisfies property (ii).}

First, to ensure property (i) holds, we require at each $\zeta\in X_0$ that
\begin{equation}\label{eq:homogeneous}
a(H_{k,\zeta})_z+b(H_{k,\zeta})_w\not\equiv0,
\end{equation}
where $H_{k,\zeta}$ denotes the first nonzero homogeneous term in the Taylor expansion of
$f(z,w)$ centered at $(\zeta,\bar\zeta)$, where $k\ge2$.
These conditions (one condition for each point in $X_0$) exclude finitely many proper
complex linear subspaces, leaving a full-measure set of parameters.

At $\zeta\in X_0$, the homogeneous term of
degree $k-1$ in $f_{\pm\xi}$ is
\[
\pm\bigl(a(H_{k,\zeta})_z+b(H_{k,\zeta})_w\bigr),
\]
which is nonzero by \eqref{eq:homogeneous}. Thus, both maps have
order of vanishing exactly $k-1$ at $\zeta$. This shows that property (i) is satisfied.

As a second requirement on $\xi = (a,b)$, at each $\zeta\in X_1$, we impose
\[
a f_z(\zeta,\bar\zeta)+b f_w(\zeta,\bar\zeta)
\notin \operatorname{im}_{\R}D F_f(\zeta),
\qquad F_f(z)=f(z,\bar z).
\]
The map $(a,b)\mapsto
a f_z(\zeta,\bar\zeta)+b f_w(\zeta,\bar\zeta)$
is onto $\C$, whereas the image of $D F_f(\zeta)$ is
a real line. Hence, these conditions also exclude only
finitely many proper real linear subspaces, leaving a full-measure set of parameters.

Choose pairwise disjoint neighborhoods $U_\zeta$ of the points
$\zeta\in X_1$, avoiding $X_0\cup X_2$.
Fix a direction $\eta=(a_0,b_0)$ satisfying both sets of
linear conditions above, and put $G=a_0f_z+b_0f_w$.
Applying Lemma~\ref{lem:rankone} in each $U_\zeta$, we obtain,
for all sufficiently small $t>0$, the sum of the number of nondegenerate zeros of $f_{t\eta}$ and
$f_{-t\eta}$ in each $U_\zeta$ is at least two.
Summing over the neighborhoods $U_\zeta$ gives at least $2|X_1|$
nondegenerate zeros. Hence, one 
of $f_{t\eta}$ and
$f_{-t\eta}$ contributes at least $|X_1|$ nondegenerate zeros, proving property~(ii).

Property (iii) holds for any sufficiently small perturbation, since nondegenerate zeros are stable.

The zeros just obtained establishing (ii) and (iii) are nondegenerate, so they persist
for parameters $\xi$ in a sufficiently small neighborhood of $t\eta$.   Shrinking the neighborhood if necessary also preserves property (i) regarding the order reductions
at the points of $X_0$. 
It remains to choose $\xi$ in this neighborhood so that property (iv) is satisfied, i.e., the zero set is finite, and every zero outside $X_0$ is nondegenerate.
We accomplish this using parametric transversality,
showing that the latter property holds on a full measure set of parameters.

At a zero of $f_\xi(z,\bar z)$ outside $X_0$, at least one
of $f_z(z,\bar z)$ and $f_w(z,\bar z)$ is nonzero.
Since the map
\[
(z,\xi)\longmapsto
f(z,\bar z)+a f_z(z,\bar z)+b f_w(z,\bar z)
\]
depends (affine) linearly on the parameters $\xi=(a,b)$, its
differential in the $\xi$-variables is the linear map
\[
(\alpha,\beta)\longmapsto
\alpha f_z(z,\bar z)+\beta f_w(z,\bar z),
\]
which is onto $\C$ at every such zero.
Proposition \ref{prop:parametric-transversality}
therefore implies that, except for a measure zero set of parameters
$\xi\in\C^2$, all zeros of $f_\xi$ outside $X_0$ are
nondegenerate, verifying property (iv). Intersecting this full measure set with its
negative (which also has full measure) imposes the same condition on $f_{-\xi}$.
For these parameters, the zero set of $f_{\xi}$ (as well as $f_{-\xi}$) is finite: the zeros outside the (finite) set $X_0$ are isolated, so they are connected components of the corresponding semialgebraic
zero set, which has only finitely many connected components
\cite[Theorem 2.4.5]{BCR98}.

Choosing $\xi$ in the intersection of this full-measure set with the neighborhood constructed above proves the claim,
and hence establishes the lemma.
\end{proof}

With the above lemma in hand, we may regularize the zeros of $f(z,\bar z)$; once they are regularized, the zeros persist under further perturbation.  Hence, it remains only to establish that the conditions stated in Lemma \ref{lem:generic} hold for an open dense set of parameters.  This is done in the next lemma, where the main step in the proof is to show that a set of ``bad'' cases forms a proper algebraic subvariety of the parameter space (following a similar structure to the discriminant argument in \cite[Theorem 3.1]{ACLRP21}).

\begin{lem}\label{lem:generic}
For $n,m\ge1$, an open dense set of polynomials with
$\deg_z f\le n$ and $\deg_w f\le m$ have bihomogenizations
whose zero curves $C\subset \CP^1 \times \CP^1$ are smooth,
avoid $(\infty,\infty)$, and are transverse to both $L$
and $\tau C$.
\end{lem}

\begin{proof}

Let $V_{n,m}$ be the space of bihomogeneous polynomials
\[
 P(Z,W)=\sum_{i=0}^n\sum_{j=0}^m
 c_{ij}Z_0^{n-i}Z_1^iW_0^{m-j}W_1^j
\]
of bidegree $(n,m)$. This family has no base points, since the four monomials
$Z_0^nW_0^m$, $Z_1^nW_0^m$, $Z_0^nW_1^m$, and
$Z_1^nW_1^m$ have no common zero in $\CP^1\times\CP^1$.
Hence, Bertini's theorem \cite[Ch.~III, Cor.~10.9]{Har77}
implies that $C$ is smooth for an open dense set of coefficients. 

The condition $c_{nm}\ne0$, which also holds on an open dense set,
ensures that $C$ avoids $(\infty,\infty)$.

For transversality of $C$ and $\tau C$, we will consider the pair of polynomials defining $C$ and $\tau C$ as a point
 $(f,g) \in V_{n,m} \times V_{m,n}$.
Non-transversality of $C$ and $\tau C$ is equivalent to the vanishing of three polynomials
in the variables $z,w$ and $c_{ij}$:
\[
f, g, \text{   and } \left| \begin{array}{cc} f_z & f_w \\ g_z & g_w\end{array}\right|.
\]
Projecting this ``non-transverse locus'' from $\CP^1\times \CP^1 \times V_{n,m} \times V_{m,n}$ onto $V_{n,m} 
\times V_{m,n}$ defines a proper closed algebraic subvariety $X$. While $X$ is a complex variety, we can
consider it as an algebraic subvariety of $\R^{4(m+1)(n+1)}$ to observe its interaction with $\tau$.
The set of points 
$$L_\tau=\{(f,g): f^* = g\} \subseteq V_{n,m} \times V_{m,n}$$ is a $2(m+1)(n+1)$-dimensional real-linear subspace of $\R^{4(m+1)(n+1)}$. To show $L_\tau$ is not contained in $X$, it is sufficient to exhibit a polynomial $f$ for which $C,\tau C$ are in fact transverse.

The example we give is a grid of horizontal and vertical lines: \[
f(z,w) = \prod_{i = 1}^n (z-a_i) \prod_{i=1}^m (w-b_i) \implies f^*(z,w) = \prod_{i = 1}^n (w-\bar a_i) \prod_{i=1}^m (z-\bar b_i).
\] 
As long as $a_i$ and $b_i$ are chosen so that, together with their conjugates, they form a set of $2(m+n)$ distinct complex numbers, the 
grids defined by $f$ and $f^*$ are transverse. This example demonstrates that for any bidegree $(n,m)$, the subspace $L_\tau \not\subseteq X$.
This shows that $X\cap L_\tau$ is a proper real algebraic subset of
$L_\tau$.  As such, its complement in $L_\tau$ is open and dense.  Thus, the generic $f$ defines a curve $C$ that is transverse to its reflection $\tau C$

Finally, $C$ is also transverse to $L$. Indeed, if
$x\in C\cap L$ and $v\in T_xC\cap T_xL$, then
$d\tau_x(v)=v$, so
\[
v\in d\tau_x(T_xC)=T_x(\tau C).
\]
Then transversality of $C$ and $\tau C$ gives $v=0$. Thus, all the required conditions hold
on an open dense set of coefficients.
\end{proof}

Lemma~\ref{lem:regularization} replaces the original polynomial by one
with at least as many zeros, all nondegenerate. A further sufficiently
small perturbation satisfying Lemma~\ref{lem:generic} preserves these
zeros and their nondegeneracy while imposing the desired additional genericity conditions stated in Lemma \ref{lem:generic-reduction}. As explained above, applying the bounds proved in Section~\ref{sec:generic}
to this generic polynomial proves
Theorem~\ref{thm:linear} in full.

\section{Examples showing asymptotic sharpness}\label{sec:examples}

The following examples, whose zero sets consist of intersections between a family of lines and a family of concentric circles (cf. \cite{Oka18}), show that, for every fixed $m$, Theorem~\ref{thm:linear} is asymptotically optimal as $n \rightarrow \infty$.

\begin{prop}\label{prop:radial-examples}
For each fixed $m\ge1$ and all sufficiently large $n$, there is a
polynomial $f(z,w)$ of bidegree $(n,m)$ for which $f(z,\bar z)$ has
finitely many zeros and $N(f)\ge(4m+1)(n-m)$.
\end{prop}

\begin{proof}
Put $k=n-m$ and consider $F(z) := f(z,\bar z)$ of the form 
\[
 F(z) =z^kA(|z|^2)+B(|z|^2),
\]
where $A,B$ are real polynomials of degree $m$ to be chosen.
Writing $z=re^{i\theta}$ and taking the imaginary part of the equation $F(r e^{i\theta})=0$ gives
\[
 r^kA(r^2)\sin(k\theta)=0,
\]
which consists of $k$ lines through the origin
and circles where $A(r^2)=0$. As $A,B$ will be chosen to have no common roots,
these circles contain no zeros of $F$.  Hence, every zero of $F$ lies on one of the alternating rays where $r>0$ and $e^{ik\theta}=\pm1$ (the origin is not a zero of $F$ since we will choose $B$ with $B(0) \neq 0$). Restricting $F(z)$ to these rays, we have 
\[
 F(re^{i\pi q/ k}) = \begin{cases} B(r^2) + r^kA(r^2), & q \text{  even}. \\ 
 B(r^2) - r^kA(r^2), & q \text{  odd}. \end{cases}
\]
Let $H_+(r) := B(r^2) + r^kA(r^2)$ denote the case where $q$ is even, and $H_-(r) := B(r^2) - r^kA(r^2)$ the $q$ odd case. Each positive root of either radial polynomial $H_\pm$ gives $k$ zeros of $F$.

Choose $B$ monic with $m$ simple roots in $(0,1)$, and $A$ with
$m$ simple roots in $(1,\infty)$, normalized so that $A(0)>0$.
Fix such a pair of polynomials while letting $k\to\infty$. Then
\begin{equation}\label{eq:Hcvgc}
 H_\pm(r)\longrightarrow B(r^2)\quad(0<r<1),
 \qquad
 r^{-k}H_\pm(r)\longrightarrow\pm A(r^2)\quad(r>1),
\end{equation}
uniformly on compact subintervals of the indicated regions.
Choose small disjoint intervals around the roots of $B(r^2)$ and
$A(r^2)$, with opposite signs at their endpoints.
For sufficiently large $k$, the intermediate value theorem gives
at least $m$ roots of each $H_\pm$ below $1$ and another $m$ above $1$.

There is an additional ``middle'' root for $H_-$ situated between these two groups of zeros. Choose $u<1<v$
with $u$ beyond all the inner intervals and $v$ before all the outer
intervals. We have $B(u^2)>0$ since $B$ is monic, and we have $A(v^2)>0$ by the normalization chosen for $A$.  Hence, the same convergence \eqref{eq:Hcvgc} gives
$H_-(u)>0>H_-(v)$ for sufficiently large $k$.
Hence, the intermediate value theorem supplies an additional root in $(u,v)$.

We therefore obtain at least
$k(4m+1)=(4m+1)(n-m)$ distinct zeros of $F$. The zero set is finite, since the zeros 
arise from the roots of nonzero univariate polynomials $H_\pm$.
Finally, $F(z)=f(z,\bar z)$ for
$f(z,w)=z^kA(zw)+B(zw)$, whose bidegree is $(k+m,m)=(n,m)$.
\end{proof}

\bibliographystyle{amsalpha}
\bibliography{polyanalytic}
\end{document}